\documentclass[]{interact}

\usepackage{epstopdf}
\usepackage[caption=false]{subfig}

\usepackage[numbers,sort&compress]{natbib}
\usepackage{graphicx}%
\usepackage{multirow}%
\usepackage{amsmath,amssymb,amsfonts}%
\usepackage{amsthm}%
\usepackage{mathrsfs}%
\usepackage[title]{appendix}%
\usepackage{xcolor}%
\usepackage{textcomp}%
\usepackage{manyfoot}%
\usepackage{booktabs}%
\usepackage{algorithm}%
\usepackage{algorithmicx}%
\usepackage{algpseudocode}%
\usepackage{listings}%
\usepackage{lineno}
\usepackage{txfonts} 
\usepackage{multirow} 
\usepackage{caption} 
\usepackage{txfonts} 
\usepackage{array} 
\usepackage{booktabs} 
\usepackage{longtable}
\usepackage{bm}
\usepackage{setspace} 
\usepackage{natbib} 
\usepackage{hyperref}
\usepackage{diagbox}
\usepackage{makecell}
\usepackage{longtable}

\usepackage{algorithm}    
\usepackage{algpseudocode}
\makeatletter
\def\BState{\State\hskip-\ALG@thistlm}
\makeatother

\bibpunct[, ]{[}{]}{,}{n}{,}{,}
\makeatletter
\def\NAT@def@citea{\def\@citea{\NAT@separator}}
\makeatother

\theoremstyle{plain}
\newtheorem{theorem}{Theorem}[section]

\theoremstyle{definition}

\theoremstyle{remark}
\newtheorem{remark}{Remark}

\begin{document}


\title{A Spectral Projected Gradient Method for Computational Protein Design problem}

\author{
\name{Yukai Zheng\textsuperscript{a} and Qingna Li\textsuperscript{b} \thanks{CONTACT Qingna Li. Email: qnl@bit.edu.cn}}
\affil{\textsuperscript{a}School of Mathematics and Statistics, Beijing Institute of Technology, Beijing, China; \textsuperscript{b}School of Mathematics and Statistics / Beijing Key Laboratory on MCAACI, Beijing Institute of Technology, Beijing, China}
}

\maketitle

\begin{abstract}
In this paper, we consider the computational protein design (CPD) problem, which is usually modeled as 0/1 programming and is extremely challenging due to its combinatorial properties. As a quadratic semi-assignment problem (QSAP), the CPD problem has been proved to be equivalent to its continuous relaxation problem (RQSAP), in terms of sharing the same optimal objective value. However, since the current algorithm for solving this RQSAP uses the projected Newton method, which requires direct computation of the Hessian matrix, its computational cost remains quite high. Precisely for this reason, we choose to employ the spectral projected gradient (SPG) method to solve the CPD problem, whose effectiveness relies on choosing the step lengths according to novel ideas that are related to the spectrum of the underlying local Hessian. Specifically, we apply the SPG method in two distinct ways: direct solving the relaxation problem and applying a penalty method. Numerical results on benchmark instances verify the superior performance of our approach over the current algorithms in both quality and efficiency.
\end{abstract}


\section{Introduction}\label{sec1}

The computational protein design (CPD) problem arises from biology, which attempts to guide the protein design process by producing a set of specific proteins that is not only rich in functional proteins, but also small enough to be evaluated experimentally. In this way, the problem of selecting amino acid sequences to perform a given task can be defined as a computable optimization problem. It is often described as the inverse of the protein folding problem \cite{pabo1983molecular,chiu1998optimizing,yue1992inverse}: the three-dimensional structure of a protein is known, and we need to find the amino acid sequence folded into it \cite{creighton1990protein}.

The challenge of the CPD problem lies in its combinatorial properties over different choices of natural amino acids. The resulting optimization model is usually NP-hard \cite{pierce2002protein,traore2013new}. Existing methods for CPD problems make use of different mathematical models, including probabilistic graphical model \cite{thomas2008protein,gainza2016algorithms}, integer linear programming model \cite{zhu2007mixed,lippow2007progress}, 0/1 quadratic programming model \cite{riazanov2017inverse,forrester2008quadratic} and weighted partial maximum satisfiability problem (MaxSAT) \cite{luo2017ccehc,schiex2014computational}. Various models were proposed in different situations with different scopes. However, due to the exponential complexity, these branch-and-bound approaches cannot solve  large-scale CPD problems.
Therefore, some preprocessing methods were proposed to reduce the problem size and improve the solution efficiency \cite{shah2004preprocessing,allouche2014computational,yanover2007dead}. For example, the dead end elimination (DEE) method \cite{allouche2014computational,yanover2007dead} reduces the problem size by eliminating some selection choices in the combinatorial space which does not contain the optimal solution. Such strategies can speed up the algorithm when sovling the CPD problems \cite{allouche2014computational}, but the worst-case complexity of the algorithm itself has not decreased.  

In recent work \cite{CPD}, it has been proved that as a quadratic semi-assignment problem (QSAP), the CPD problem is equivalent to its continuous relaxation problem (RQSAP), in terms of sharing the same optimal objective value.  Furthermore, one can obtain the global optimizer of the CPD problem by a global optimizer of its relaxation. As a consequence, compared with the existing branch-and-bound approaches that suffer from high computational complexity, the proposed algorithm AQPPG in \cite{CPD} is based on a continuous problem and enjoys a lower per-iteration complexity, which makes it suitable for solving large-scale CPD problems. This method of transforming the original problem into its equivalent continuous relaxation problem is also applicable to any assignment model where the objective function linearly depends on its variables, such as hypergraph matching \cite{hypergraph}, mimo detection \cite{mimo} and graph clustering problems \cite{clustering}. However, since the proposed algorithm AQPPG uses the projected Newton method to solve the subproblem, which requires direct computation of the Hessian matrix of the objective function, its computational cost remains quite high, indicating potentials for improvement. Meanwhile, we notice that the spectral projected gradient (SPG) method \cite{birgin2000nonmonotone, birgin2001algorithm, birgin2003inexact}, whose effectiveness relies on choosing the step lengths according to novel ideas that are related to the spectrum of the underlying local Hessian, can effectively utilize the second-order information of the objective function, while avoid direct computation of the Hessian matrix. This motivates our work in this paper. The contributions of this paper is summarized as follows. Firstly, we choose to use the SPG method to solve the subproblem of the CPD problem. Secondly, we apply the SPG method in two distinct ways, namly SCSC and SCP. SCSC directly applies the spectral projected gradient method to solve the relaxation of the CPD problem, while SCP employs a penalty method to solve the penalized relaxation.

The rest of this paper is organized as follows. In Section 2, we introduce the CPD problem and its equivalent relaxation. In Section 3, we introduce the idea of spectral projected gradient method and propose two distinct approaches, SCSC and SCP, to solve the CPD problem using the SPG method. In Section 4, we report the numerical results. Final conclusions are made in Section 5.

\section{The CPD Problem}\label{sec2}

In this section, we introduce the CPD problem and its equivalent relaxation.

\subsection{Problem Formulation}

\begin{figure}
\centering
\includegraphics[width=0.6\linewidth]{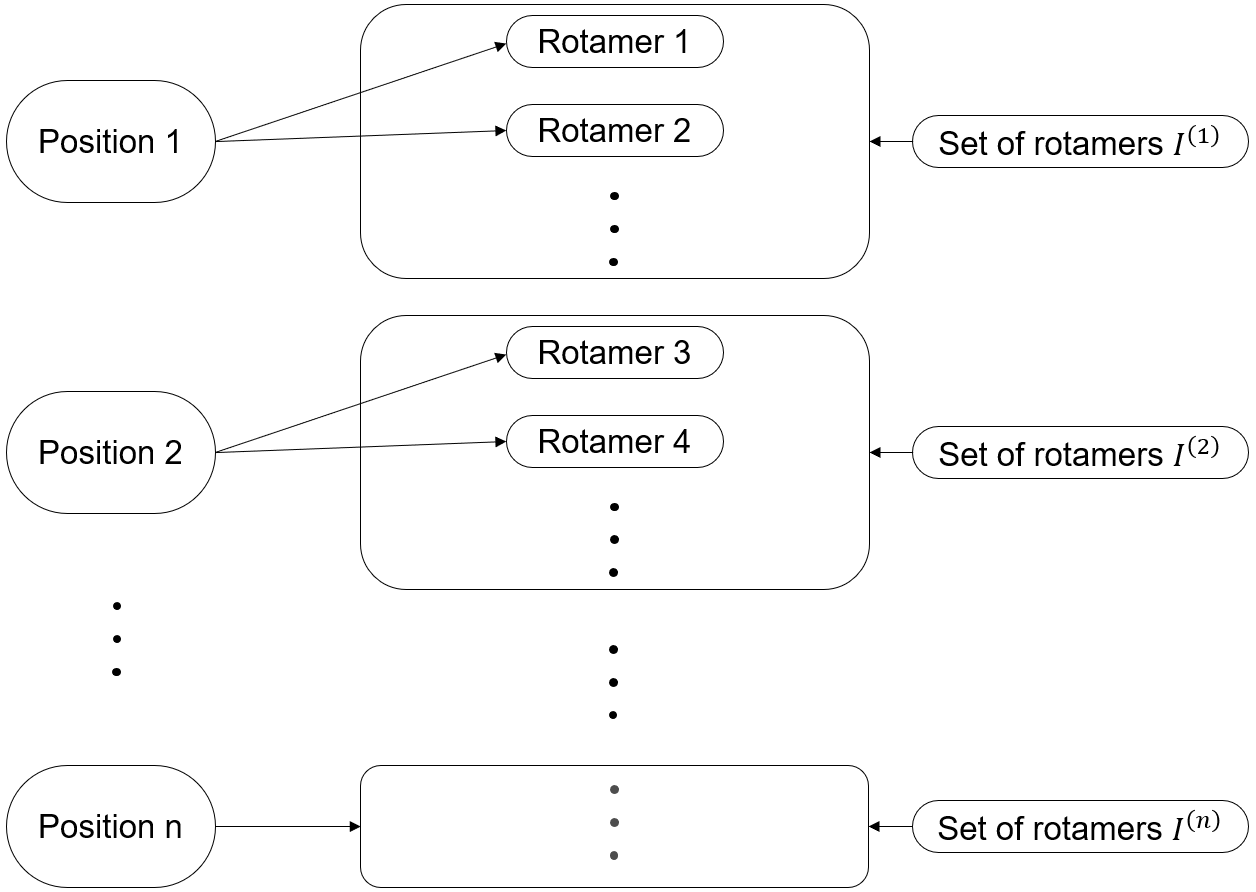}
\caption{The relationship between positions and rotamers in the CPD problem}\label{p_r}
\end{figure}

Briefly, the CPD problem is to select a rotamer among a set of them for each position on the backbone to minimize the total energy, as shown in \textbf{Figure} \ref{p_r}. To be more specific, the CPD problem could be described as the following mathematical model.

Let $n$ be the number of positions on the backbone of the protein, and we use $[n]$ to denote the set of integers $\left\{1,2,...,n\right\}$. $I^{(i)}$ is a set of rotamers that one can choose for position $i\in [n]$, where $i_r$ corresponds to the $r$-th rotamer that can be chosen for position $i$. Let $l_i=\vert I^{(i)}\vert$, that is, the number of elements in the set $I^{(i)}$. Let $m=\sum\limits^n_{i=1}l_i$. Define $x\in\mathbb{R}^m$ as follows:

\begin{equation}\label{vectorx}
x=
\begin{pmatrix}
x^{(1)}\\ \vdots\\x^{(n)}
\end{pmatrix}
=
\begin{pmatrix}
x^{(1)}_1\\ \vdots \\x^{(1)}_{l_1}\\ \vdots\\x^{(n)}_1\\ \vdots \\x^{(n)}_{l_n}
\end{pmatrix}
\in\mathbb{R}^m.\notag
\end{equation}
Here $x^{(i)}\in\mathbb{R}^{l_i}$ is the $i$-th block of the assignment variable $x$, $i\in [n]$, which satisfies
\begin{equation}
x^{(i)}_r=
\begin{cases}
1,\enspace\hbox{if the }r\hbox{-th rotamer is assigned to position }i,\enspace i\in [n],\\
0,\enspace\hbox{otherwise.}\notag
\end{cases}
\end{equation}
The CPD problem is to pick up one specific rotamer for each position of the protein, such that the total energy of the protein is minimized, that is,

\begin{equation}\label{assignment}
\begin{split}
\mathop{\min}_{x\in\mathbb{R}^m}\enspace& f(x):=\frac{1}{2}x^\top Bx+a^\top x,\\
\hbox{s.t.}\enspace&\sum\limits_{r\in [l_i]}x^{(i)}_r=1,\enspace i\in [n],\\
&x^{(i)}_r\in \left\{0,1\right\},\enspace r\in [l_i],\enspace i\in [n].
\end{split}
\end{equation}
where
\begin{equation}\label{vectora}
a=
\begin{pmatrix}
a^{(1)}\\ \vdots\\a^{(n)}
\end{pmatrix}
=
\begin{pmatrix}
a^{(1)}_1\\ \vdots \\a^{(1)}_{l_1}\\ \vdots\\a^{(n)}_1\\ \vdots \\a^{(n)}_{l_n}
\end{pmatrix}
\in\mathbb{R}^m,\enspace B=
\begin{pmatrix}
\bf{0}_{l_1\times l_1} & B_{12} & \cdots & B_{1n} \\
B_{12}^\top  & \bf{0}_{l_2\times l_2} & \cdots & B_{2n} \\
\vdots & \vdots & \ddots & \vdots \\
B_{1n}^\top  & B_{2n}^\top  & \cdots & \bf{0}_{l_n\times l_n} \\
\end{pmatrix}
\in\mathbb{R}^{m\times m},\notag
\end{equation}
\begin{equation}
B_{ij}=
\begin{pmatrix}
b^{ij}_{11} & b^{ij}_{12} & \cdots & b^{ij}_{1l_j}\\
b^{ij}_{21} & b^{ij}_{22} & \cdots & b^{ij}_{2l_j}\\
\vdots & \vdots & \ddots & \vdots\\
b^{ij}_{l_i1} & b^{ij}_{l_i2} & \cdots & b^{ij}_{l_il_j}\\
\end{pmatrix}
\in\mathbb{R}^{l_i\times l_j},\enspace b^{ij}_{rs}=E(i_r,j_s),\enspace i<j,\enspace i,\enspace j\in [n].\notag
\end{equation}
$a_i^r$ represents the energy contribution of rotamer $r$ at position $i$ capturing internal interactions (and a reference energy for the associated amino acid) or interactions with fixed regions, and $b_{rs}^{ij}$ represents the pairwise interaction energy between rotamer $r$ at position $i$ and rotamer $s$ at position $j$ \cite{desmet1992dead}. $a_i^r$ and $b_{rs}^{ij}$ are constant and already known for $i$, $j\in [n]$, $r\in I^{(i)}$, $s\in I^{(j)}$. 

\subsection{An Equivalent Continuous Relaxation Problem}

Consider the following relaxation of the CPD problem (\ref{assignment}):
\begin{equation}\label{relaxation}
\begin{split}
\mathop{\min}_{x\in\mathbb{R}^{m}}\enspace&f(x)\\
\hbox{s.t.}\enspace&\sum\limits_{r\in [l_i]}x^{(i)}_r=1,\enspace i\in [n],\\
&x^{(i)}_r\geq 0,\enspace r\in [l_i],\enspace i\in [n].
\end{split}
\end{equation}
It has been proved that the CPD problem (\ref{assignment}) is equivalent to this continuous relaxation problem (\ref{relaxation}), in terms of sharing the same optimal objective value \cite{CPD}.  Specifically, we have the following result. 

\begin{theorem}\cite[Theorem 3.1]{CPD}\label{theorem1}
Let $\Vert x\Vert_0$ denote the number of nonzero elements in $x$.
There exists an optimal solution 
\begin{math}
x^*
\end{math}
to the relaxation problem (\ref{relaxation}) such that 
\begin{math}
\Vert x^*\Vert_0=n.
\end{math}
Moreover,
\begin{math}
x^*
\end{math}
is also an optimal solution to the original problem {\rm (\ref{assignment})}.
\end{theorem}

\textbf{Theorem} \ref{theorem1} basically reveals that problem (\ref{relaxation}) is a tight continuous relaxation of problem (\ref{assignment}) in the sense that two problems share at least one global minimizer.  Furthermore, given a global minimizer of the relaxation problem (\ref{relaxation}), one can use the following algorithm to obtain the global minimizer of the CPD problem (\ref{assignment}). 

\begin{algorithm}\label{algorithm1}
\caption{Obtain the global minimizer of (\ref{assignment}) by that of (\ref{relaxation})}\label{algorithm1}
\begin{algorithmic}[1]
\Require a global optimal solution $x^0\in\mathbb{R}^m$ of the relaxation problem (\ref{relaxation}).
\Ensure a global optimal solution $\hat{x}\in\mathbb{R}^m$ of the original problem (\ref{assignment}).
\State initialization: $l=0$.
\While{$\Vert x^l\Vert_0>n$}

		\State for $i=1,...,n$, find the first block of $x^l$, denoted as ${(x^l)}^{(j)}$, such that $\Vert {(x^l)}^{(j)}\Vert_0>1$, choose one index $s^0$ from $\Gamma_j(x^l):=\left\{s:{(x^l)}^{(j)}_s>0\right\}$, and define $x^{l+1}$ as:

		\begin{math}
		(x^{l+1})^{(j)}_s=
		\begin{cases}
		1, &s=s^0,\\
		0, &\hbox{otherwise,}
		\end{cases}
		\enspace {(x^{l+1})}^{(i)}=
		\begin{cases}
		{(x^{l+1})}^{(j)}, &i=j,\\
		{(x^l)}^{(i)}, &\hbox{otherwise.}
		\end{cases}
		\end{math}
        
		\State $l=l+1$.
\EndWhile\\
output $\hat{x}=x^l.$
\end{algorithmic}
\end{algorithm}

\begin{remark}
For any approximate solution of the relaxation problem (\ref{relaxation}), we can also use Algorithm \ref{algorithm1} to obtain a feasible solution of the original problem (\ref{assignment}). 
\end{remark}

The convergence of Algorithm \ref{algorithm1} is provided in \cite[Theorem 4.3]{CPD}. Based on the results above, in the rest of this paper, we obtain the solution of the original CPD problem (\ref{assignment}) by simply solving its equivalent continuous relaxation (\ref{relaxation}).

\section{The SPG Method and Two Approaches for (\ref{relaxation})}\label{sec3}

In this section,  we introduce the idea of spectral projected gradient method and propose two distinct approaches SCSC and SCP to solve the CPD problem using the SPG method.

\subsection{The Framework of SPG Method}

The spectral projected gradient (SPG) method \cite{birgin2000nonmonotone, birgin2001algorithm, birgin2003inexact} for solving convex constrained problems was born from the marriage of the global Barzila-Borwein (spectral) nonmonotone scheme with the classical projected gradient (PG) method, which have been extensively used in statistics. Indeed, the effectiveness of the classical PG method can be greatly improved by incorporating the spectral step length and nonmonotone globalization strategies. 

Now we introduce the specific process of the SPG method \cite{birgin2014spectral}. Consider the following convex constrained problem,

\begin{equation}\label{convex}
\begin{split}
\mathop{\min}_{x}\enspace& g(x)\\
\hbox{s.t.}\enspace& x\in\Omega,
\end{split}
\end{equation}
where $g$ is a convex function and $\Omega$ is a closed convex set. Let $P_\Omega$ denote the Euclidean projection onto the set $\Omega$. The rough idea for SPG is to update $x^{k+1}$ at iteration $k$ by $x^{k+1}=x^K+\alpha_kd^k$, where $d^k=P_\Omega(x^k-\nabla g(x^k))-x^k$ and $d^k$ is given by nonmonotone line search, whose initial trial is the BB stepsize. Assume that a stopping criterion parameter $\epsilon>0$, a sufficient decrease parameter $\gamma\in (0, 1)$, an integer parameter $M\geq 1$ for the nonmonotone line search, safeguarding parameters $0<\sigma_1<\sigma_2<1$ for the quadratic interpolation, safeguarding parameters $0<\lambda_{min}\leq\lambda_{max}<\infty$ for the spectral step length, an arbitrary initial $x^*$, and $\lambda_0\in [\lambda_{min}, \lambda_{max}]$ are given. Details of SPG method are described in Algorithm \ref{spg}.

\begin{algorithm}\label{spg}
\caption{The spectral projected gradient method}\label{spg}
\begin{algorithmic}[1]
\Require $x^*$, $P_\Omega$, $\epsilon$, $\gamma$, $M$, $\sigma$, $\sigma_2$, $\lambda_{min}$, $\lambda_{max}$, $\lambda_0$.
\Ensure a global optimal solution $\overline{x}\in\Omega$ of the convex constrained problem (\ref{convex}).
\State initialization: if $x^*\notin\Omega$, redefine $x^* = P_\Omega(x^*)$, $k := 0$.
\While{$\Vert P_\Omega(x^k-\nabla g(x^k))-x^k\Vert_\infty>\epsilon$}

		\State set the search direction $d_k=P_\Omega(x^k-\nabla g(x^k))-x^k$, compute the step length $\alpha_k$ using Algorithm \ref{linesearch} (with parameters $\gamma$, $M$, $\sigma_1$ and $\sigma_2$), and set $x^{k+1}=x^k+\alpha_kd_k$.

		\State compute $s_k=x^{k+1}-x^k$ and $t_k=\nabla g(x^{k+1})-\nabla g(x^k)$. If $s_k^\top t_k\leq0$, then set $\lambda_{k+1}=\lambda_{max}$. Otherwise, set $\lambda_{k+1}=max\left\{\lambda_{min}, min\left\{s_k^\top s_k/s_k^\top t_k, \lambda_{max}\right\}\right\}$.

		$k=k+1$.
\EndWhile\\
output $\overline{x}=x^k.$
\end{algorithmic}
\end{algorithm}

\begin{algorithm}\label{linesearch}
\caption{Nonmonotone line search}\label{linesearch}
\begin{algorithmic}[1]
\Require $x^k$, $d_k$, $\gamma$, $M$, $\sigma_1$, $\sigma_2$, $g(x^{(k-j)})$, $\forall j\in [min\left\{k, M-1\right\}]$.
\Ensure step length $\alpha_k$.
\State initialization: compute $g_{max}=max\left\{g(x^{k-j})\ |\ 0\leq j\leq min\left\{k, M-1\right\}\right\}$, and set $\alpha = 1$.
\While{$g(x^k+\alpha d_k)> g_{max}+\gamma\alpha\nabla g(x^k)^\top d_k$}

		\State compute $\alpha_{tmp}=-\frac{1}{2}\alpha^2\nabla g(x^k)^\top d_k/[g(x^k+\alpha d_k)-g(x^k)-\alpha\nabla g(x^k)^\top d_k]$. If $\alpha_{tmp}\in [\sigma_1, \sigma_2\alpha]$, then set $\alpha=\alpha_{tmp}$. Otherwise, set $\alpha=\alpha/2$.		

\EndWhile\\

output $\alpha_k=\alpha.$
\end{algorithmic}
\end{algorithm}

\subsection{Two Approaches Based on the SPG Method}

\subsubsection{SCSC}

The SPG method is suitable for solving the CPD problem in this paper, for it can effectively utilize the second-order information of the objective function, while avoid direct computation of the Hessian matrix.
As a consequence, we propose two distinct approaches to solve the CPD problem using the SPG method. The first one, namely SCSC, is to solve the relaxation problem (\ref{relaxation}) directly and then obtain the solution of the original CPD problem (\ref{assignment}). Note that the constraints in (\ref{relaxation}) are typical simplex constraints \cite{boyd2004convex}, and the projection $P_\Omega$ is not that obvious. We use the method in \cite{chen2011projection}, as shown in Algorithm \ref{projection}, to obtain the projection of any points in $\mathbb{R}^{m}$ onto the feasible region
\begin{equation}\label{sets}
\Omega_s:=\left\{x\ |\ \sum\limits_{r}x^{(i)}_r=1, x^{(i)}_r\geq 0,\forall r\in [l_i],\forall i\in [n]\right\},
\end{equation}
whose computational complexity is $O(nlog(n))$. The algorithm determines which components should remain positive after projection by sorting the components in order of magnitude, then sets the remaining components to zero. During this process, it reasonably sets the threshold $\hat{t}$ to ensure that the sum of the truncated components is exactly $1$. The essence of Algorithm \ref{projection} is to evenly distribute the mass that exceeds the feasible region to the active components. The specific process of the whole approach can be represented as Algorithm \ref{SCSC} below. We term this approch as SCSC, which is the abbreviation of SPG for CPD with Simplex Constraints.
\begin{algorithm}\label{projection}
\caption{Compute the projection of $x$ onto the feasible region $\Omega_s$}\label{projection}
\begin{algorithmic}[1]
\Require $x=(x^{(1)},...,x^{(n)})^\top\in\mathbb{R}^m$.
\Ensure the projection $P_s(x)$.
\For {$k$ = 1, 2,..., n}

		\State sort $x^{(k)}$ in the ascending order as $x^{(k)}_1\leq...\leq x^{(k)}_{l_k}$, and set $i=l_k-1$.

		\State compute $t_i=\frac{\sum^{l_k}_{j=i+1}x^{(k)}_j-1}{l_k-i}$. If $t_i\geq x^{(k)}_i$ then set $\hat{t}=t_i$ and go to Step 5, otherwise set $i=i-1$ and redo Step 3 if $i\geq 1$ or go to Step 4 if $i=0$.

		\State set $\hat{t}=\frac{\sum^{l_k}_{j=1}x^{(k)}_j-1}{l_k}$.

		\State set $P_s(x^{(k)})=(x^{(k)}-\hat{t})_+$.

\EndFor\\

output $P_s(x)=(P_s(x^{(1)}),...,P_s(x^{(n)}))^\top.$
\end{algorithmic}
\end{algorithm}
\begin{algorithm}\label{SCSC}
\caption{SCSC}\label{SCSC}
\begin{algorithmic}[1]
\Require $x^*$, $P_\Omega$, $\epsilon$, $\gamma$, $M$, $\sigma$, $\sigma_2$, $\lambda_{min}$, $\lambda_{max}$, $\lambda_0$.
\Ensure a global optimal solution $\hat{x}$ of the CPD problem (\ref{assignment}).
\State initialization: if $x^*\notin\Omega_s$, redefine $x^* = P_{\Omega_s}(x^*)$, $k := 0$.\\

		use Algorithm \ref{spg} to get the solution $\overline{x}$ of the relaxation problem (\ref{relaxation}), where $g(x)=f(x)$, $\Omega_s$ is defined as in (\ref{sets}), and $P_s(x)$ is computed by Algorithm \ref{projection}.\\

		obtain the global optimal solution $\hat{x}$ of the CPD problem (\ref{assignment}) by $\overline{x}$ using Algorithm \ref{algorithm1}.\\

output $\hat{x}.$
\end{algorithmic}
\end{algorithm}

\subsubsection{SCP}

On the other hand, it should be noticed that the aim we solve (\ref{relaxation}) is to identify the locations of nonzero entries of the global minimizer of (\ref{relaxation}), rather than to find the magnitude of it. This is because once the locations of the nonzero entries are identified, we can apply Algorithm \ref{algorithm1} to obtain a global optimal solution of (\ref{assignment}). From this point of view, keeping the equality constraints in (\ref{relaxation}) may not be necessary. Therefore, we apply the quadratic penalty method to solve (\ref{relaxation}). That is, we penalize the equality constraints to the objective function, and solve the following quadratic penalty subproblem,
\begin{equation}\label{penalty}
\begin{split}
\mathop{\min}_{x\in\mathbb{R}^m}\enspace& f(x)+\frac{\sigma}{2}\sum\limits_{i=1}^n\left(\sum\limits_{r\in [l_i]}x^{(i)}_r-1\right)^2\\
\hbox{s.t. }\enspace& x\in\Omega_P,
\end{split}
\end{equation}
where $\sigma>0$ is a penalty parameter, and
\begin{equation}\label{setp}
\Omega_P =\left\{x\ |\ x^{(i)}_r\geq 0,\forall r\in [l_i],\forall i\in [n]\right\}.
\end{equation} 
As a consequence, anothor way to solve the CPD problem is to solve the penalty problem (\ref{penalty}) and then obtain the solution of the original problem (\ref{assignment}) as shown in Algorithm \ref{SCP}, where $x_+$ denotes the positive part of $x$. We term this approch as SCP, which is the abbreviation of SPG for CPD Penalized.
\begin{algorithm}\label{SCP}
\caption{SCP}\label{SCP}
\begin{algorithmic}[1]
\Require $x^*$, $P_\Omega$, $\epsilon$, $\gamma$, $M$, $\sigma$, $\sigma_2$, $\lambda_{min}$, $\lambda_{max}$, $\lambda_0$.
\Ensure a global optimal solution $\hat{x}$ of the CPD problem (\ref{assignment}).
\State initialization: if $x^*\notin\Omega_p$, redefine $x^* = P_{\Omega_p}(x^*)$, $k := 0$.\\

		use Algorithm \ref{spg} to get the solution $\overline{x}$ of the penalty problem (\ref{penalty}), where $g(x)=f(x)+\frac{\sigma}{2}\sum\limits_{i=1}^n\left(\sum\limits_{r\in [l_i]}x^{(i)}_r-1\right)^2$, $\Omega_P$ is defined as in (\ref{setp}), and $P_p(x) = x_+$.\\

		obtain the global optimal solution $\hat{x}$ of the CPD problem (\ref{assignment}) by $\overline{x}$ using Algorithm \ref{algorithm1}.\\

output $\hat{x}.$
\end{algorithmic}
\end{algorithm}

Compared with the existing algorithm AQPPG \cite{CPD} for solving the CPD problem through its relaxation, which uses the projected Newton method and requires direct computation of the Hessian matrix of the objective function, our algorithms, both SCSC and SCP, could effectively utilize the second-order information of the objective function, while avoid direct computation of the Hessian matrix. Specifically, since the step sizes of our algorithms are determined by the spectrum of the local Hessian matrix, it can be expected that the efficiency of our algorithms, particularly in terms of the number of iterations, will be significantly better than that of the existing algorithm \cite{CPD}.

\section{Numerical Results}\label{sec5}

We implement the algorithm in MATLAB (R2024b). All the experiments are performed on a ROG Asus desktop with AMD Ryzen 9 7940HX CPU at 2.4 GHz and 16 GB of memory running Windows 11. We use the data as in \cite{allouche2014computational}, which can be downloaded from \url{https://genoweb.toulouse.inra.fr/~tschiex/CPD-AIJ/}.
\footnote{To convert the floating point energies of a given instance to non-negative integer costs, David Allouche et al. \cite{allouche2014computational} subtracted the minimum energy to all energies and then multiplied energies by an integer constant $Q$ and rounded to the nearest integer. Therefore, all the energies in the data sets are non-negative integers.}
\textbf{Table} \ref{parameter} shows the information of all data sets tested.  In \textbf{Table} \ref{parameter}, $n$ represents the number of positions in the target protein, which is also the number of blocks in the decision variable $x$. $\min l_i$ represents how many components one block in $x$ contains at least, and $\max l_i$ represents how many components one block contains at most. $m$ shows the dimension of the decision variable $x$.

\begin{longtable}[H]{|c|c|c|c|c|c|}
\caption{Information of all the data sets}\label{parameter}
\endfirsthead
\multicolumn{6}{c}{{Table \ref{parameter} - continued from previous page}}\\ \hline
NO. & Data & $n$ & $l_i\in (\min l_i, \max l_i)$ & $m$\\\hline
\endhead
\hline\multicolumn{3}{l}{{Continued on next page}}\\
\endfoot
\hline
\endlastfoot
\hline
NO. & Data & $n$ & $l_i\in (\min l_i, \max l_i)$ & $m$\\\hline
1 & 1HZ5 & 12 & (49, 49) & 588 \\ 
2 & 1PGB & 11 & (49, 49) & 539 \\ 
3 & 2PCY & 18 & (48, 48) & 864 \\ 
4 & 1CSK & 30 & (3, 49) & 616 \\ 
5 & 1CTF & 39 & (3, 56) & 1204 \\ 
6 & 1FNA & 38 & (3, 48) & 990 \\ 
7 & 1PGB & 11 & (198, 198) & 2178 \\ 
8 & 1UBI & 13 & (49, 49) & 637 \\ 
9 & 2TRX & 11 & (48, 48) & 528 \\ 
10 & 1UBI & 13 & (198, 198) & 2574 \\ 
11 & 2DHC & 14 & (198, 198) & 2772 \\ 
12 & 1PIN & 28 & (198, 198) & 5544 \\ 
13 & 1C9O & 55 & (198, 198) & 10890 \\ 
14 & 1C9O & 43 & (3, 182) & 1950 \\ 
15 & 1CSE & 97 & (3, 183) & 1355 \\ 
16 & 1CSP & 30 & (3, 182) & 1114 \\ 
17 & 1DKT & 46 & (3, 190) & 2243 \\ 
18 & 1BK2 & 24 & (3, 182) & 1294 \\ 
19 & 1BRS & 44 & (3, 194) & 3741 \\ 
20 & 1CM1 & 17 & (198, 198) & 3366 \\ 
21 & 1SHG & 28 & (3, 182) & 737 \\ 
22 & 1MJC & 28 & (3, 182) & 493 \\ 
23 & 1SHF & 30 & (3, 56) & 638 \\ 
24 & 1FYN & 23 & (3, 186) & 2474 \\ 
25 & 1NXB & 34 & (3, 56) & 800 \\ 
26 & 1TEN & 39 & (3, 66) & 808 \\ 
27 & 1POH & 46 & (3, 182) & 943 \\ 
28 & 1CDL & 40 & (3, 186) & 4141 \\ 
29 & 1HZ5 & 12 & (198, 198) & 2376 \\ 
30 & 2DRI & 37 & (3, 186) & 2120 \\ 
31 & 2PCY & 46 & (3, 56) & 1057 \\ 
32 & 2TRX & 61 & (3, 186) & 1589 \\ 
33 & 1CM1 & 42 & (3, 186) & 3633 \\ 
34 & 1LZ1 & 59 & (3, 57) & 1467 \\ 
35 & 1GVP & 52 & (3, 182) & 3826 \\ 
36 & 1R1S & 56 & (3, 182) & 3276 \\ 
37 & 2RN2 & 69 & (3, 66) & 1667 \\ 
38 & 1HNG & 85 & (3, 182) & 2341 \\ 
39 & 3CHY & 74 & (3, 66) & 2010 \\ 
40 & 1L63 & 83 & (3, 182) & 2392 \\ 
\end{longtable}

\subsection{Details of the experimental setup}

Here we address two details of settings of the numerical experiments. The first one is about the termination conditions. It should be emphasised that the aim we solve (\ref{relaxation}) is to identify the locations of nonzero entries of the global minimizer of (\ref{relaxation}), rather than to find the magnitude of it. This is because once the locations of the nonzero entries are identified, we can apply Algorithm \ref{algorithm1} to obtain a global optimal solution of (\ref{assignment}). From this point of view, we need to formulate reasonable termination criteria, so that the algorithm stops at the appropriate time, rather than iterating to complete convergence. In practice, we stop the iteration once any of the following holds:
\begin{itemize}
\item [(a)] $\vert f(x^{k+1})-f(x^k)\vert<\epsilon_a$,
\item [(b)] $\max\limits_{j\in \left\{0,1,...,M-1\right\}}\vert f(x^{k+j})-f(x^{k+j-M})\vert<\epsilon_b$,
\item [(c)] For every iteration, namely $k$, we use Algorithm \ref{algorithm1} to transform $x^k$ into a feasible point $\hat{x^k}$ of the original CPD problem (\ref{assignment}). We stop the iteration if $\hat{x^k}$ stays the same for $N$ iterations.
\end{itemize}
Criterion (a) is a commonly-used terminal condition. As for (b), we take $M$ iterations as a unit to observe the behavior of the variable $x$ in order to avoid falling into a dead loop, since the SPG method imposes a functional decrease every $M$ iterations. The last criterion (c) aims to stop the iterations if the locations of nonzero entries of the global minimizer are already identified. It is expected that most cases are terminated by the third criterion (c). In our experiments, we set $\epsilon_a=10^{-2}$, $\epsilon_b=10^{-2}$, $M=10$ and $N=50$.

The second detail is about the nonmonotone line search, which we use to solve the subproblem of the SPG method. It is suggested to set the initial value of $\alpha$ to 1 in Nonmonotone line search. However, we find that setting $\alpha$ to different values, such as $0.9$, $1.1$ and $2$, may leads to better solutions. We hypothesize that this is because setting $\alpha$ to $1$ introduces some cases that are too special. For example, the variable $x$ may turn to exactly zero after the first iteration, making the initial point $x^0$ invalid. As a consequence, we set the initial value of $\alpha$ to $0.9$ in our experiments.

\subsection{Comparison with the state-of-the-art algorithm}

\begin{figure}[H]
\centering
\includegraphics[width=1\linewidth]{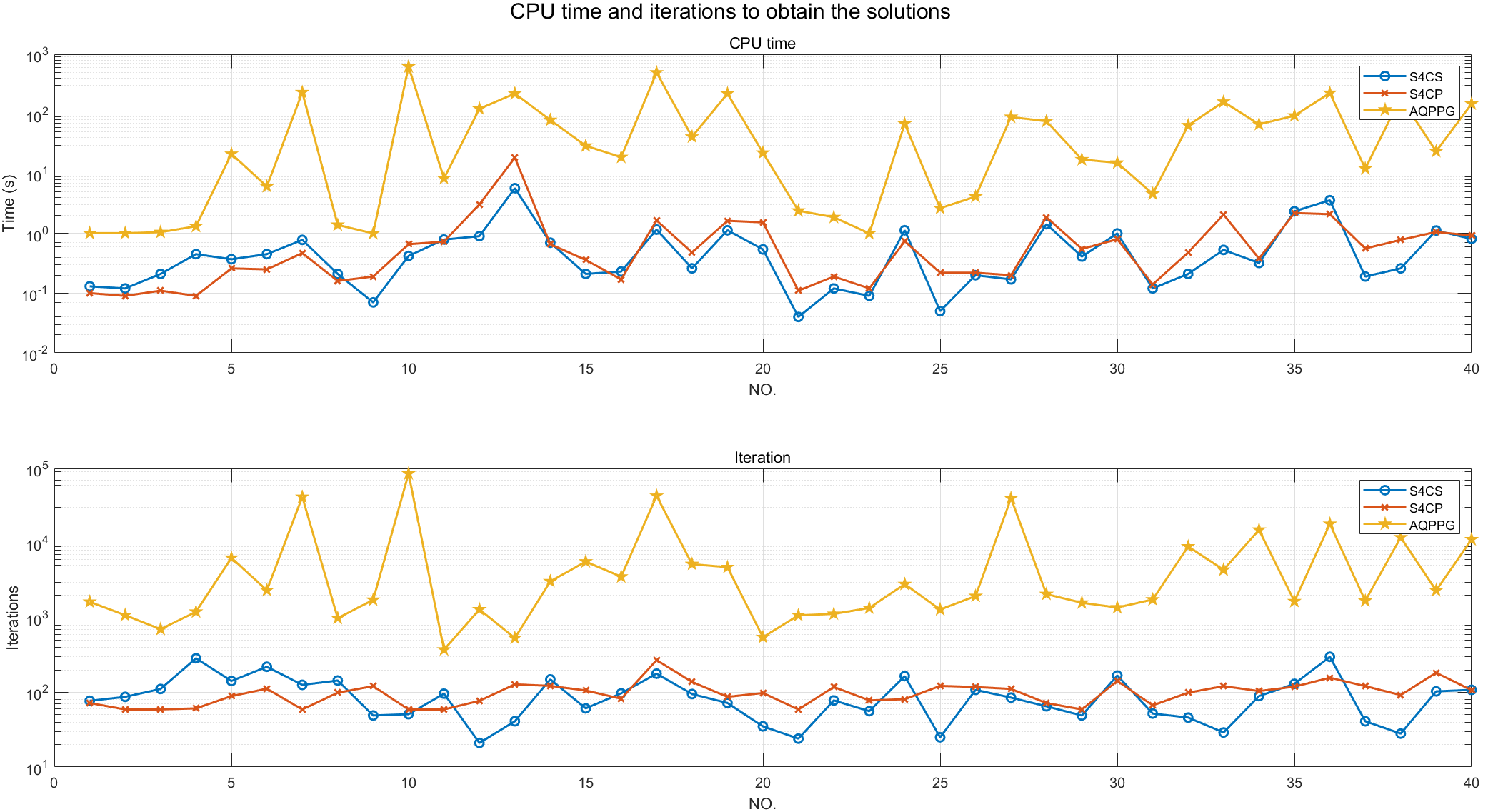}
\caption{CPU time and iterations for each method}\label{efficiency}
\end{figure}

We compare SCSC and SCP with AQPPG, which is the state-of-the-art algorithm \cite{CPD} for solving the CPD problem through its continuous relaxation (\ref{relaxation}). 
\textbf{Table} \ref{results} shows the results given by AQPPG and Gurobi.
$m$ shows the dimension of the decision variable $x$.
$Objective$ represents the function values given by different methods. 
$Ratio$ represents the ratio of the optimal values given by SCSC and SCP compared to those given by AQPPG. 
\textbf{Figure} \ref{efficiency} illustrates the efficiency of different algorithms. The upper and lower graphs respectively show the CPU time spent and the number of iterations during the solving process. In \textbf{Figure} \ref{efficiency}, the dotted line with circles represents SCSC, the dotted line with crosses represents SCP, and the dotted line with five-pointed stars represents AQPPG. \textbf{Table} \ref{time} shows detailed data of \textbf{Figure} \ref{efficiency}. $Iteration$ and $Time$ respectively represent the number of iterations and the CPU time spent to obtain the solutions. The unit of time is seconds (accurate to two decimal places).
\begin{longtable}[H]{|c|c|c|rrr|rr|}
\caption{Results given by SCSC, SCP and AQPPG}\label{results}
\endfirsthead
\multicolumn{6}{c}{{Table \ref{results} - continued from previous page}}\\ \hline
\multirow {2} {0.75cm} {NO.} & \multirow {2} {0.75cm} {Data} & \multirow {2} {0.3cm} {m} & \multicolumn {3} {c|} {Objective} & \multicolumn {2} {c|} {Ratio} \\
&&& \multicolumn {1} {c} {SCSC} & \multicolumn {1} {c} {SCP} & \multicolumn {1} {c|} {AQPPG} & \multicolumn {1} {c} {SCSC} & \multicolumn {1} {c|} {SCP} \\ \hline
\endhead
\hline\multicolumn{3}{l}{{Continued on next page}}\\
\endfoot
\hline
\endlastfoot
\hline
\multirow {2} {0.75cm} {NO.} & \multirow {2} {0.75cm} {Data} & \multirow {2} {0.3cm} {m} & \multicolumn {3} {c|} {Objective} & \multicolumn {2} {c|} {Ratio} \\
&&& \multicolumn {1} {c} {SCSC} & \multicolumn {1} {c} {SCP} & \multicolumn {1} {c|} {AQPPG} & \multicolumn {1} {c} {SCSC} & \multicolumn {1} {c|} {SCP} \\ \hline
1 & 1HZ5 & 588 & 150978 & 153821& 151818 & 99.45\% & 101.32\%  \\ 
2 & 1PGB & 539 & 125509 & 127854 & 125509 & 100.00\% & 101.87\% \\  
3 & 2PCY & 864 & 308094 & 308691 & 308545 & 99.85\% & 100.05\% \\  
4 & 1CSK & 616 & 1125970 & 1132771 & 1126030 & 99.99\% & 100.60\% \\ 
5 & 1CTF & 1204 & 1882233 & 1889393 & 1883085 & 99.95\% & 100.33\% \\ 
6 & 1FNA & 990 & 3752105 & 3758744 & 3751852 & 100.01\% & 100.18\% \\ 
7 & 1PGB & 2178 & 286140 & 288601 & 288170 & 99.30\% & 100.15\% \\  
8 & 1UBI & 637 & 160053 & 163033 & 160053 & 100.00\% & 101.86\% \\ 
9 & 2TRX & 528 & 178774 & 180154 & 178705 & 100.04\% & 100.81\% \\  
10 & 1UBI & 2574 & 381117 & 384126 & 385034 & 98.98\% & 99.76\% \\ 
11 & 2DHC & 2772 & 1422426 & 1425825 & 1443284 & 98.55\% & 98.79\% \\ 
12 & 1PIN & 5544 & 1996302 & 1998834 & 1996834 & 99.97\% & 100.10\% \\ 
13 & 1C9O & 10890 & 8018607 & 8023657 & 8084802 & 99.18\% & 99.24\% \\ 
14 & 1C9O & 1950 & 4962718 & 4989300 & 4975017 & 99.75\% & 100.29\% \\ 
15 & 1CSE & 1355 & 18603227 & 18652539 & 18646035 & 99.77\% & 100.03\% \\ 
16 & 1CSP & 1114 & 2523524 & 2544597 & 2529459 & 99.77\% & 100.60\% \\ 
17 & 1DKT & 2243 & 4194065 & 4237681 & 4214282 & 99.52\% & 100.56\% \\ 
18 & 1BK2 & 1294 & 1138719 & 1143602 & 1143207 & 99.61\% & 100.03\% \\ 
19 & 1BRS & 3741 & 4008271 & 4030471 & 4028208 & 99.51\% &100.06\% \\ 
20 & 1CM1 & 3366 & 743672 & 746489 & 746221 & 99.66\% & 100.04\% \\ 
21 & 1SHG & 737 & 1513822 & 1515142 & 1514734 & 99.94\% & 100.03\% \\ 
22 & 1MJC & 493 & 1514559 & 1517143 & 1518774 & 99.72\% & 99.89\% \\ 
23 & 1SHF & 638 & 1101866 & 1104045 & 1101912 & 99.99\% & 100.19\% \\ 
24 & 1FYN & 2474 & 1183942 & 1191313 & 1194046 & 99.15\% & 99.77\% \\ 
25 & 1NXB & 800 & 2978188 & 3005473 & 2978209 & 99.99\% & 100.92\% \\ 
26 & 1TEN & 808 & 1960254 & 1963554 & 1960322 & 99.99\% & 100.16\% \\ 
27 & 1POH & 943 & 4034124 & 4045543 & 4034259 & 99.99\% & 100.28\% \\ 
28 & 1CDL & 4141 & 3592304 & 3668267 & 3620297 & 99.23\% & 101.33\% \\ 
29 & 1HZ5 & 2376 & 343863 & 345548 & 345034 & 99.66\% & 100.15\% \\ 
30 & 2DRI & 2120 & 2906028 & 2916207 & 2956196 & 98.30\% & 98.64\% \\ 
31 & 2PCY & 1057 & 2937093 & 2949181 & 2937992 & 99.97\% & 100.38\% \\ 
32 & 2TRX & 1589 & 7017563 & 7053178 & 7022829 & 99.93\% & 100.43\% \\ 
33 & 1CM1 & 3633 & 3897482 & 3913161 & 3904719 & 99.81\% & 100.22\% \\ 
34 & 1LZ1 & 1467 & 7027973 & 7048132 & 7031221 & 99.95\% & 100.24\% \\ 
35 & 1GVP & 3826 & 5198362 & 5214654 & 5205320 & 99.87\% & 100.18\% \\ 
36 & 1R1S & 3276 & 6174978 & 6228964 & 6177881 & 99.95\% & 100.83\% \\ 
37 & 2RN2 & 1667 & 8918391 & 8970316 & 8922415 & 99.95\% & 100.54\% \\ 
38 & 1HNG & 2341 & 13539761 & 13590905 & 13543984 & 99.97\% & 100.35\% \\ 
39 & 3CHY & 2010 & 10466525 & 10473661 & 10568323 & 99.04\% & 99.10\% \\ 
40 & 1L63 & 2392 & 12892684 & 12960382 & 13015089 & 99.06\% & 99.58\% \\ 
\end{longtable}
\begin{longtable}[H]{|c|c|c|rrr|rrr|}
\caption{CPU time and iterations to obtain the solutions}\label{time}
\endfirsthead
\multicolumn{9}{c}{{Table \ref{time} - continued from previous page}}\\ \hline
\multirow {2} {0.75cm} {NO.} & \multirow {2} {0.75cm} {Data} & \multirow {2} {0.3cm} {$m$} & \multicolumn {3} {c|} {Time} & \multicolumn {3} {c|} {Iteration} \\
&& \multicolumn {1} {c|} {} & \multicolumn {1} {c} {SCSC} & \multicolumn {1} {c} {SCP} & \multicolumn {1} {c|} {AQPPG} & \multicolumn {1} {c} {SCSC} & \multicolumn {1} {c} {SCP} & \multicolumn {1} {c|} {AQPPG} \\ \hline
\endhead
\hline\multicolumn{9}{l}{{Continued on next page}}\\
\endfoot
\hline
\endlastfoot
\hline
\multirow {2} {0.75cm} {NO.} & \multirow {2} {0.75cm} {Data} & \multirow {2} {0.3cm} {$m$} & \multicolumn {3} {c|} {Time} & \multicolumn {3} {c|} {Iteration} \\
&& \multicolumn {1} {c|} {} & \multicolumn {1} {c} {SCSC} & \multicolumn {1} {c} {SCP} & \multicolumn {1} {c|} {AQPPG} & \multicolumn {1} {c} {SCSC} & \multicolumn {1} {c} {SCP} & \multicolumn {1} {c|} {AQPPG} \\ \hline

1 & 1HZ5 & 588 & 0.13 & 0.10 & 1.01 & 77 & 72 & 1631  \\ 
2 & 1PGB & 539 & 0.12 & 0.09 & 1.01 & 87 & 59 & 1081 \\  
3 & 2PCY & 864 & 0.21 & 0.11 & 1.05 & 111 & 59 & 701 \\  
4 & 1CSK & 616 & 0.45 & 0.09 & 1.32 & 286 & 61 & 1201 \\ 
5 & 1CTF & 1204 & 0.37 & 0.26 & 21.49 & 142 & 89 & 6337 \\ 
6 & 1FNA & 990 & 0.45 & 0.25 & 6.12 & 220 & 112 & 2326 \\ 
7 & 1PGB & 2178 & 0.78 & 0.47 & 230.63 & 126 & 59 & 41493 \\  
8 & 1UBI & 637 & 0.21 & 0.16 & 1.39 & 144 & 99 & 988 \\ 
9 & 2TRX & 528 & 0.07 & 0.19 & 1.00 & 49 & 121 & 1731 \\  
10 & 1UBI & 2574 & 0.42 & 0.66 & 621.83 & 51 & 59 & 85220 \\ 
11 & 2DHC & 2772 & 0.79 & 0.73 & 8.38 & 96 & 59 & 373 \\ 
12 & 1PIN & 5544 & 0.90 & 3.07 & 122.50 & 21 & 77 & 1291 \\ 
13 & 1C9O & 10890 & 5.73 & 18.53 & 220.08 & 41 & 128 & 534 \\ 
14 & 1C9O & 1950 & 0.70 & 0.66 & 79.20 & 149 & 122 & 3068 \\ 
15 & 1CSE & 1355 & 0.21 & 0.36 & 29.35 & 61 & 106 & 5639 \\ 
16 & 1CSP & 1114 & 0.23 & 0.17 & 18.93 & 97 & 82 & 3541 \\ 
17 & 1DKT & 2243 & 1.16 & 1.65 & 493.61 & 178 & 270 & 42841 \\ 
18 & 1BK2 & 1294 & 0.26 & 0.48 & 41.58 & 95 & 138 & 5217 \\ 
19 & 1BRS & 3741 & 1.13 & 1.63 & 219.54 & 72 & 87 & 4749 \\ 
20 & 1CM1 & 3366 & 0.54 & 1.53 & 22.40 & 35 & 98 & 548 \\ 
21 & 1SHG & 737 & 0.04 & 0.11 & 2.40 & 24 & 59 & 1076 \\ 
22 & 1MJC & 493 & 0.12 & 0.19 & 1.87 & 78 & 118 & 1124 \\ 
23 & 1SHF & 638 & 0.09 & 0.12 & 1.00 & 56 & 78 & 1351 \\ 
24 & 1FYN & 2474 & 1.13 & 0.74 & 68.54 & 165 & 81 & 2813 \\ 
25 & 1NXB & 800 & 0.05 & 0.22 & 2.65 & 25 & 122 & 1287 \\ 
26 & 1TEN & 808 & 0.20 & 0.22 & 4.14 & 108 & 118 & 1948 \\ 
27 & 1POH & 943 & 0.17 & 0.20 & 89.10 & 85 & 111 & 39738 \\ 
28 & 1CDL & 4141 & 1.43 & 1.87 & 75.73 & 65 & 72 & 2068 \\ 
29 & 1HZ5 & 2376 & 0.41 & 0.55 & 17.31 & 49 & 59 & 1584 \\ 
30 & 2DRI & 2120 & 1.00 & 0.80 & 15.22 & 168 & 142 & 1369 \\ 
31 & 2PCY & 1057 & 0.12 & 0.14 & 4.61 & 52 & 67 & 1752 \\ 
32 & 2TRX & 1589 & 0.21 & 0.48 & 64.47 & 46 & 99 & 8990 \\ 
33 & 1CM1 & 3633 & 0.53 & 2.07 & 160.41 & 29 & 121 & 4379 \\ 
34 & 1LZ1 & 1467 & 0.32 & 0.38 & 67.15 & 89 & 104 & 15029 \\ 
35 & 1GVP & 3826 & 2.35 & 2.20 & 93.75 & 130 & 119 & 1657 \\ 
36 & 1R1S & 3276 & 3.60 & 2.10 & 224.71 & 299 & 157 & 18047 \\ 
37 & 2RN2 & 1667 & 0.19 & 0.56 & 12.16 & 41 & 121 & 1684 \\ 
38 & 1HNG & 2341 & 0.26 & 0.78 & 187.57 & 28 & 92 & 11898 \\ 
39 & 3CHY & 2010 & 1.12 & 1.06 & 23.70 & 103 & 183 & 2314 \\ 
40 & 1L63 & 2392 & 0.81 & 0.92 & 148.16 & 108 & 107 & 11190 \\ 
\end{longtable}
\noindent We can see that our algorithm AQPPG could effectively solve CPD problems. In terms of efficiency, the proposed SCSC and SCP both significantly outperform the existing algorithm AQPPG. SCSC and SCP are able to solve most instances within 1 second, whereas AQPPG exceeds 60 seconds for many cases. Additionally, the numbers of iterations for SCSC and SCP are one to two orders of magnitude smaller than those of AQPPG. From the perspective of quality, the solutions obtained by SCSC generally surpass those of AQPPG, ranging from 98.30\% to 100.04\%. The solution quality of SCP is similar to that of AQPPG, with most differences being less than 1\%. Specifically, we observe that the performance of SCSC is superior to that of SCP. We speculate that this is due to the following reasons:
\begin{itemize}
\item [(1)] Due to the use of the penalty method, SCP requires determining the size of the penalty parameter $\sigma$. In our experiments, we set $\sigma = 10^7$, which is approximately equal to or slightly greater than the optimal value of the objective function. However, this may not necessarily be the optimal strategy. Setting an empirically derived optimal $\sigma$ for each instance could potentially lead to better performance.
\item [(2)] In SCP, the projection function $P_{\Omega_p}$ simply takes the maximum between the projected element itself and zero. This projection is overly simplistic and may result in a loss of significant information during the iterative search.
\item [(3)] From the perspective of termination criteria, compared to SCP, SCSC more frequently concludes iterations based on the third termination criterion (c) mentioned earlier, indicating its better adaptability to the termination criteria.
\end{itemize}

Based on the above results, we can conclude that our proposed methods, particularly SCSC, could effectively find high-quality solutions for the CPD problems within a reasonable amount of time.

\section{Conclusion}\label{sec6}

In this paper, we introduce two novel approaches SCSC and SCP for solving the CPD problem based on the spectral projected gradient method by addressing their equivalent relaxation problems. Specifically, SCSC directly applies the spectral projected gradient method to solve the relaxation of the CPD problem, while SCP employs a penalty method to solve the penalized relaxation. Benefiting from the property of the spectral projected gradient method that allows for the effective exploitation of second-order information of the objective function without directly computing the Hessian matrix, our proposed methods are particularly well-suited for solving large-scale CPD problems. Numerical results verify the superior performance of our approaches, particularly SCSC, over the current algorithms in both quality and efficiency.

\section*{Data availability statement}
The datasets used in this paper can be downloaded from \url{https://genoweb.toulouse.inra.fr/~tschiex/CPD-AIJ/}.

\section*{Disclosure statement}
The authors have no conflict of interest to declare that are relevant to the content of this article.

\section*{Funding}
The work of Qing-Na Li is supported by the National Natural Science Foundation of China (NSFC) 12071032 and 12271526.

\end{document}